\documentclass{article}
\usepackage{amsfonts}
\usepackage{amssymb}
\usepackage{amsthm}
\usepackage{latexsym}
\usepackage{xspace}
\usepackage{amsmath}
\usepackage[latin1]{inputenc}
\usepackage{epsfig}

\newtheorem{dfn}{Definition}[section]
\newtheorem{ex}{Example}[section]
\newtheorem{thm}{Theorem}[section]

\newcommand{\bb}{\boldsymbol{b}}
\newcommand{\bx}{\boldsymbol{x}}

\newcommand{\by}{\boldsymbol{y}}
\newcommand{\bz}{\boldsymbol{z}}

\newcommand{\bc}{\boldsymbol{c}}

\newcommand{\bSigma}{\boldsymbol{\Sigma}}
\newcommand{\be}{\boldsymbol{e}}
\newcommand{\bzero}{\boldsymbol{0}}

\begin{document}

\title{Array Variate Skew Normal Random Variables with Multiway Kronecker Delta Covariance Matrix Structure}
\author{Deniz Akdemir \\ Department of Statistics \\ 
  University of Central Florida\\ Orlando, FL 32816}

\maketitle

\begin{abstract} In this paper, we will discuss the concept of an array variate random variable and introduce a class of skew normal array densities that are obtained through a selection model that uses the array variate normal density as the kernel and the cumulative distribution of the univariate normal  distribution as the selection function. \end{abstract}


\allowdisplaybreaks

\section{Introduction}

The array variate random variable up to $2$ dimensions has been studied intensively in \cite{gupta2000matrix} and by many others. For arrays observations of $3,$ $4$ or in general $i$ dimensions a suitable normal probability model has been recently proposed in \cite{DenizGupta}. An elliptical generalization of the normal array variable can be found in [Akdemir(2011)].  In this paper, we will generalize the notion of elliptical random variables to the array case.

In Section 2, we first study the algebra of arrays, and introduce the concept of an array variable random variable.  In Sections 3, the density of the normal array random variable is provided. We finally provide the definition for an array random variable with skew normal density in Section 4.

\section{Array Algebra}

In this paper we will only study arrays with real elements. We will write $\widetilde{X}$ to say that $\widetilde{X}$ is an array. When it is necessary we can write the dimensions of the array as subindices, e.g., if $\widetilde{X}$ is a $m_1 \times m_2\times m_3 \times m_4$ dimensional array in $R^{m_1\times m_2 \times \ldots \times m_i}$, then we can write $\widetilde{X}_{m_1 \times m_2\times m_3 \times m_4}.$  Arrays with the usual element wise summation and scalar multiplication operations can be shown to be a vector space.  

To refer to an element of an array $\widetilde{X}_{m_1 \times m_2\times m_3 \times m_4},$ we write the position of the element as a subindex to the array name in parenthesis, $(\widetilde{X})_{r_1r_2r_3r_4}.$  If we want to refer to a specific column vector obtained by keeping all but an indicated dimension constant, we indicate the constant dimensions as before but we will put '$:$' for the non constant dimension, e.g., for  $\widetilde{X}_{m_1 \times m_2\times m_3 \times m_4},$    $(\widetilde{X})_{r_1r_2:r_4}$ refers to the the column vector $(({X})_{r_1r_21r_4},({X})_{r_1r_22r_4}, \ldots, ({X})_{r_1r_2m_3r_4})'.$
 
We will now review some basic principles and techniques of array algebra. These results and their proofs can be found in Rauhala \cite{rauhala1974array},  \cite{rauhala1980introduction}  and Blaha \cite{blaha1977few}. 

\begin{dfn}\emph{Inverse Kronecker product} of two matrices $A$ and $B$ of dimensions $p\times q$ and $r \times s$ correspondingly is written as $A\otimes^i B$ and is defined as $A\otimes^i B=[A(B)_{jk}]_{pr\times qs}=B\otimes A,$ where $'\otimes'$ represents the ordinary Kronecker product.\end{dfn}

The following properties of the inverse Kronecker product are useful:\begin{itemize}\item $\bzero \otimes^i A= A \otimes^i \bzero=\bzero.$ \item $(A_1+A_2)\otimes^i B=A_1\otimes^i B+ A_2 \otimes^i B.$ \item $A \otimes^i (B_1+ B_2)=A \otimes^i B_1+ A \otimes^i B_2.$ \item  $\alpha A \otimes^i \beta B= \alpha \beta A \otimes^i B.$ \item $(A_1 \otimes^i B_1)(A_2 \otimes^i B_2)= A_1A_2 \otimes^i B_1B_2.$ \item  $(A \otimes^i B)^{-1}=(A^{-1} \otimes^i B^{-1}).$ \item  $(A \otimes^i B)^{+}=(A^{+} \otimes^i B^{+}),$ where $A^{+}$ is the Moore-Penrose inverse of $A.$ \item $(A \otimes^i B)^{-}=(A^{-} \otimes^i B^{-}),$ where $A^{-}$ is the $l$-inverse of $A$ defined as $A^{-}=(A'A)^{-1}A'.$\item If $\{\lambda_i\}$ and $\{\mu_j\}$ are the eigenvalues with the corresponding eigenvectors $\{\bx_i\}$ and $\{\by_j\}$ for matrices $A$ and $B$ respectively, then $A\otimes^i B$ has eigenvalues $\{\lambda_i\mu_j\}$ with corresponding eigenvectors $\{\bx_i\otimes^i\by_j\}.$\item Given two matrices $A_{n\times n}$ and $B_{m\times m}$ $|A\otimes^i B|=|A|^m|B|^n,$ $tr(A\otimes^i B)=tr(A)tr(B).$\item $A\otimes^i B=B\otimes A=U_1 A\otimes B U_2,$ for some permutation matrices $U_1$ and $U_2.$ 
\end{itemize}

It is well known that a matrix equation $$AXB'=C$$ can be rewritten in its monolinear form as \begin{equation}\label{eqmnf}A\otimes^i B vec(X)=vec(C).\end{equation} Furthermore, the matrix equality $$A\otimes^i B XC'=E$$ obtained by stacking equations of the form (\ref{eqmnf}) can be written in its monolinear form as $$(A\otimes^i B \otimes^i C) vec(X)=vec(E).$$ This process of stacking equations could be continued and R-matrix multiplication operation introduced by Rauhala \cite{rauhala1974array} provides a compact way of representing these equations in array form:

\begin{dfn}\emph{R-Matrix Multiplication} is defined element wise: 

 $$((A_1)^1 (A_2)^2 \ldots (A_i)^i\widetilde{X}_{m_1 \times m_2 \times \ldots \times m_i})_{q_1q_2\ldots q_i}$$ $$=\sum_{r_1=1}^{m_1}(A_1)_{q_1r_1}\sum_{r_2=1}^{m_2}(A_2)_{q_2r_2}\sum_{r_3=1}^{m_3}(A_3)_{q_3r_3}\ldots \sum_{r_i=1}^{m_i}(A_i)_{q_ir_i}(\widetilde{X})_{r_1r_2\ldots r_i}.$$ \end{dfn}

R-Matrix multiplication generalizes the matrix multiplication (array multiplication in two dimensions)to the case of $k$-dimensional arrays. The following useful properties of the R-Matrix multiplication are reviewed by Blaha \cite{blaha1977few}:
\begin{enumerate}
\item $(A)^1B=AB.$ 
\item $(A_1)^1(A_2)^2C=A_1CA'_2.$ 
\item $\widetilde{Y}=(I)^1(I)^2\ldots (I)^i \widetilde{Y}.$ 
\item \small $((A_1)^1 (A_2)^2\ldots (A_i)^i)((B_1)^1(B_2)^2\ldots (B_i)^i)\widetilde{Y}= (A_1B_1)^1(A_2B_2)^2\ldots(A_iB_i)^i\widetilde{Y}.$ \normalsize
\end{enumerate}

The operator $rvec$ describes the relationship between $\widetilde{X}_{m_1 \times m_2 \times \ldots m_i}$ and its monolinear form $\bx_{m_1m_2\ldots m_i\times 1}.$ 
\begin{dfn}\label{def:rvec} $rvec( \widetilde{X}_{m_1 \times m_2 \times \ldots m_i})=\bx_{m_1m_2\ldots m_i\times 1}$ where $\bx$ is the column vector obtained by stacking the elements of the array $\widetilde{X}$ in the order of its dimensions; i.e., $(\widetilde{X})_{j_1 j_2 \ldots j_i}=(\bx)_j$ where $j=(j_i-1)n_{i-1}n_{i-2}\ldots n_1+(j_i-2)n_{i-2}n_{i-3}\ldots n_1+\ldots+(j_2-1)n_1+j_1.$\end{dfn}

\begin{thm}\label{rmultvec}Let $\widetilde{L}_{m_1 \times m_2 \times\ldots m_i}=(A_1)^1(A_2)^2\ldots(A_i)^i\widetilde{X}$ where $(A_j)^j$ is an $m_j\times n_j$ matrix for $j=1,2,\ldots,i$ and $\widetilde{X}$ is an $n_1\times n_2\times\ldots\times n_i$ array. Write $\mathbf{l}=rvec(\widetilde{L})$ and $\bx=rvec(\widetilde{X}).$ Then, $\mathbf{l}=A_1\otimes^iA_2\otimes^i\ldots\otimes^i A_i\bx.$ \end{thm}

Therefore, there is an equivalent expression of the array equation in monolinear form.

\begin{dfn}{} The square norm of $\widetilde{X}_{m_1 \times m_2 \times\ldots m_i}$ is defined as $$\|\widetilde{X}\|^2=\sum_{j_1=1}^{m_1}\sum_{j_2=1}^{m_2}\ldots\sum_{j_i=1}^{m_i}((\widetilde{X})_{j_1j_2\ldots j_i})^2.$$ \end{dfn}

\begin{dfn}{} The distance of $\widetilde{X_1}_{m_1 \times m_2 \times\ldots m_i}$ from $\widetilde{X_2}_{m_1 \times m_2 \times\ldots m_i}$ is defined as $\sqrt{\|\widetilde{X_1}-\widetilde{X_2}\|^2}.$ \end{dfn}

\begin{ex} Let $\widetilde{Y}=(A_1)^1 (A_2)^2\ldots (A_i)^i\widetilde{X}+\widetilde{E}.$ Then $\|\widetilde{E}\|^2$ is minimized for $\widehat{\widetilde{X}}=(A_1^{-})^1(A_2^{-})^2\ldots(A_i^{-})^i\widetilde{Y}.$ \end{ex}

\section{Array Variate Normal Distribution}

By using the results in the previous section on array algebra, mainly the relationship of the arrays to their monolinear forms described by Definition \ref{def:rvec} , we can write the density of the standard normal array variable.
\begin{dfn}{} If $$\widetilde{Z}\sim N_{m_1 \times m_2\times \ldots \times m_i}(\widetilde{M}=\bzero, \Lambda=I_{m_1m_2\ldots m_i}),$$ then $\widetilde{Z}$ has array variate standard normal distribution. The pdf of $\widetilde{Z}$ is given by \begin{equation}\label{eq:densitystarn}f_{\widetilde{Z}} (\widetilde{Z})=\frac{\exp{(-\frac{1}{2}\|\widetilde{Z}\|^2)}}{(2\pi)^{m_1m_2\ldots m_i/2}}.\end{equation}\end{dfn}

For the scalar case, the density for the standard normal variable $z \in \mathbf{R}^1$ is given as \[\phi_1(z)=\frac{1}{(2\pi)^\frac{1}{2}}exp(-\frac{1}{2}z^2).\] For the $m_1$ dimensional  standard normal vector $\bz\in\mathbf{R}^{m_1},$ the density is given by \[\phi_{m_1}(\bz)=\frac{1}{(2\pi)^\frac{m_1}{2}}exp(-\frac{1}{2}\bz'\bz).\] Finally the $m_1\times m_2$ standard matrix variate variable $Z\in\mathbf{R}^{m_1\times m_2}$ has the density \[\phi_{m_1\times m_2}(Z)=\frac{1}{(2\pi)^\frac{m_1m_2}{2}}exp(-\frac{1}{2}trace(Z'Z)).\] With the above definition, we have generalized the notion of normal random variable to the array variate case.

\begin{dfn}{\label{stdnorm}} We write $$\widetilde{X}\sim N_{m_1 \times m_2\times \ldots \times m_i}(\widetilde{M}, \Lambda_{m_1m_2\ldots m_i})$$ if   $rvec(\widetilde{X})\sim N_{m_1m_2\ldots m_i}(rvec(\widetilde{M}),\Lambda_{m_1m_2\ldots m_i}).$ Here, $\widetilde{M}$ is the expected value of $\widetilde{X}$, and $\Lambda_{m_1m_2\ldots m_i}$ is the covariance matrix of the ${m_1m_2\ldots m_i}$-variate random variable $rvec(\widetilde{X}).$\end{dfn} 

The family of normal densities with Kronecker delta covariance structure are obtained by considering the densities obtained by the location-scale transformations of the standard normal variables. This kind of model is defined in the next.

\begin{dfn}\label{modkroncov}{(\cite{DenizGupta})} Let $\widetilde{Z}\sim N_{m_1 \times m_2\times \ldots \times m_i}(\widetilde{M}=\bzero, \Lambda=I_{m_1m_2\ldots m_i}).$ Define $\widetilde{X}=(A_1)^1 (A_2)^2 \ldots (A_i)^i\widetilde{Z}+\widetilde{M}$ where $A_1, A_2,\ldots,A_i$ are non singular matrices of orders $m_1, m_2,\ldots, m_i$.   Then the pdf of $\widetilde{X}$ is given by \small \begin{equation}\label{eq:densityarn}\phi(\widetilde{X}; \widetilde{M},A_1,A_2,\ldots A_i)=\frac{\exp{(-\frac{1}{2}\|{(A_1^{-1})^1 (A_2^{-1})^2 \ldots (A_i^{-1})^i(\widetilde{X}-\widetilde{M})}\|^2)}}{(2\pi)^{m_1m_2\ldots m_i/2}|A_1|^{\prod_{j\neq 1}{m_j}} |A_2|^{\prod_{j\neq 2}{m_j}} \ldots |A_i|^{\prod_{j\neq i}{m_j}}}.\end{equation} \normalsize \end{dfn}

Distributional properties of a array normal variable with density in the form of Definition \ref{modkroncov} variable can obtained by using the equivalent monolinear representation of the random variable. The moments, the marginal and conditional distributions, independence of variates can be studied considering the equivalent monolinear form of the array variable and the well known properties of the multivariate normal random variable. 

\section{Array Variate Skew Normal Variable}

A very general definition of skew symmetric variable for the matrix case can be obtained from matrix variate selection models. Suppose $X$ is a $k\times n$ random matrix with density $f(X),$ let $g(X)$ be a weight function. A weighted form of density $f(X)$ is given by \begin{equation}\label{weightedmatrix}h(X)=\frac{f(X)g(X)}{\int_{\mathbb{R}^{k\times n}}g(X)f(X)dX}.\end{equation} When the sample is only a subset of the population then the associated model would be called a selection model. 

Chen and Gupta provide a skew normal matrix variate probability density function in the following form  (\cite{guptachen2}):\begin{equation}f_1(X,\bSigma\otimes\Psi,\bb)=c_1^{*}\phi_{k\times n}(X;\bSigma\otimes\Psi)\Phi_n(X'\bb,\Psi)\end{equation} where $c_1^{*}=(\Phi_n(\bzero, (1+\bb'\bSigma\bb)\Psi))^{-1}$, $\phi(.)$  and $\Phi(.)$ denote the density function the cumulative distribution of the standard normal random variable correspondingly. A drawback of this definition is that it allows independence only over its rows or columns, but not both. Harrar (\cite{harrar}) give two more matrix variate skew normal densities: \begin{equation}\label{eqharrar1}f_2(X,\bSigma,\Psi,\bb,\Omega)=c_2^{*}\phi_{k\times n}(X;\bSigma\otimes\Psi)\Phi_n(X'\bb,\Omega)\end{equation} and \begin{equation}\label{eqharrar2}f_3(X,\bSigma,\Psi,\bb,B)=c_3^{*}\phi_{k\times n}(X;\bSigma\otimes\Psi)\Phi_n(tr(B'X))\end{equation} 
where $c_2^{*}=(\Phi_n(\bzero, (\Omega+\bb'\bSigma\bb)\Psi))^{-1},$ $c_3^{*}=2;$ $\bSigma,$ $\Psi,$ and $\Omega$ are positive definite covariance matrices of dimensions $k,$ $n$ and $n$ respectively, $B$ is a matrix of dimension $k\times n.$ Note that if $\Omega=\Psi$ then $f_2$ is the same as $f_1.$ Although, more general than $f_1,$ the densities $f_2$ and $f_3$ still do not permit independence of rows and columns simultaneously. In \cite{akdemir2009class}, the following density is given to overcome this problem:
\small
\begin{equation}\label{msn}f(X,M,A,B,{\Delta})=\frac{\phi_{k\times n}({A}^{-1}(X-M){B}^{-1})\prod_{j=1}^{k}\prod_{i=1}^{n}\Phi(\alpha_{ji}{\be}'_{j}({A}^{-1}(X-M){B}^{-1})\bc_i)}{|A|^n|B|^k2^{-kn}}.\end{equation} \normalsize

Here, the parameter $M$ is a $k\times n$ matrix. The scale parameters $A$ and $B$ are $k\times k$ and $n \times n$ positive definite matrices correspondingly. The shape parameter $\Delta$ is a matrix of order  $k\times n.$ Finally, ${\be}_{j}$ is the unit length $k$ dimensional vector with $1$ at its jth row and  ${\bc}_{i}$ is the unit length $n$ dimensional vector with $1$ at its ith row. We denote the distribution of a variable with density of the form (\ref{msn}) by $msn_{k\times n}({M}, A,B, \Delta).$
Motivated by the model in (\ref{weightedmatrix}) we first define a general class of selection models for the array variate random variables.

\begin{dfn}Let $f(.)$ be a density for a array random variable of dimensions $m_1 \times m_2 \times\ldots \times m_i$ and let $g(.)$ be a weight function that takes an array of the same dimensions as its argument.   A weighted form of density $f(.)$ is given by \begin{equation}\label{weightedmatrix}h(\widetilde{X})=\frac{f(\widetilde{X})g(\widetilde{X})}{\int_{\mathbb{R}^{m_1 \times m_2 \times\ldots \times m_i}}g(\widetilde{X})f(\widetilde{X})d\widetilde{X}}.\end{equation} When the sample is only a subset of the population then the associated model would be called a selection model. Here $f(.)$ is called the kernel density and $g(.)$ is called the weight or selection function.\end{dfn}

We will use the density for the array variate random variable in as the kernel density for the selection model. We obtain a class of skew normal array variate densities by using a selection function like the one in (\ref{msn}):

\begin{dfn}A density for a skew normal array variate random variable $\widetilde{X}_{m_1 \times m_2 \times\ldots \times m_i}$ is given by

\begin{equation}\label{asnpdf}f(\widetilde{X};\widetilde{M},A_1,A_2,\ldots A_i,{\widetilde{\Delta}})=\frac{\phi(\widetilde{X}; \widetilde{M},A_1,A_2,\ldots A_i)g(\widetilde{X}; \widetilde{M},A_1,A_2,\ldots A_i, \widetilde{\Delta})}{2^{-m_1m_2\ldots m_i}}\end{equation} where \tiny $$g(\widetilde{X}; \widetilde{M},A_1,A_2,\ldots A_i, \widetilde{\Delta})=\prod_{j_1=1}^{m_1}\prod_{j_2=1}^{m_2}\ldots\prod_{j_i=1}^{m_i}\Phi((\widetilde{\Delta})_{j_1j_2\ldots j_i}({(A_1^{-1})^1 (A_2^{-1})^2 \ldots (A_i^{-1})^i(\widetilde{X}-\widetilde{M})})_{j_1j_2\ldots j_i}),$$ \normalsize $A_1, A_2,\ldots,A_i$ are non singular matrices of orders $m_1, m_2,\ldots, m_i,$ $\widetilde{M}$ and $\widetilde{\Delta}$ are constant arrays of the same dimension as $\widetilde{X}$.\normalsize
\end{dfn}
\section{Discussion}
The models with Kronecker delta covariance structures provide a great deal of decrease in the number of parameters that has to be estimated. We will now discuss this in the context if principle component analysis.

Principal components analysis is a useful statistical technique that has found applications in fields such as face recognition and image compression, and is a common technique for finding patterns in data of high dimension. The end product of PCA is a set of new uncorrelated variables ordered in terms of their variances obtained from a linear combination of the original variables. 

\begin{dfn} For the $m_1 \times m_2 \times\ldots \times m_i$ dimensional array variate random variable $\widetilde{X},$ the principal components are defined as the principal components of the  $d=m_1m_2\ldots m_i$-dimensional random vector $rvec(\widetilde{X}).$ \end{dfn} 

The advantage in choosing a Kronecker structure is the decrease in the number of parameters. If $\{\lambda(A_r)_{r_j}\}$ are the $m_j$ eigenvalues of $A_rA'_r$ with the corresponding eigenvectors $\{(\bx_r)_{r_j}\}$ for $r=1,2,\ldots,i$ and $r_j=1,2,\ldots,m_r,$ then $(A_1\otimes^i A_2\otimes^iA_i)(A_1\otimes^i A_2\otimes^iA_i)'$ will have  eigenvalues $\{\lambda(A_1)_{r_1}\lambda(A_2)_{r_2}\ldots\lambda(A_i)_{r_i}\}$ with corresponding eigenvectors $\{(\bx_1)_{r_1}\otimes^i(\bx_2)_{r_2}\otimes^i \ldots \otimes^i(\bx_i)_{r_i} \}.$  By replacing  $A_r$ by their estimators, we can estimate the eigenvalues and eigenvectors of the covariance of $rvec(\widetilde{X})$ using this relationship. 

When applying the array variate skew normal model to real data, more parsimonious forms of the model in Equation (\ref{asnpdf}) should be considered. For example, a model where $\widetilde{\Delta}$ is diagonal should be appropriate for most cases. One could further assume that many of these diagonal elements are zero.

\bibliographystyle{plainnat}
\bibliography{arrayref}
\end{document}